\documentclass[twoside,12pt]{article}
\usepackage{amssymb,amsmath,bm,euscript}
\usepackage{graphicx,color,caption}

\usepackage[colorlinks=true]{hyperref}

\usepackage{algorithm,algpseudocode}
\usepackage{braket}

\newtheorem{proposition}{Proposition}[section]
\newtheorem{theorem}[proposition]{Theorem}

\newtheorem{definition}[proposition]{Definition}

\newtheorem{example}[proposition]{Example}

\newcommand{\qed}{\hphantom{.}\hfill $\Box$\medbreak}

\def\O{\mathcal{O}}

\def\T{\mathcal{T}}

\def\A{{\mathcal{A}}}

\def\B{{\mathcal{B}}}
\def\C{\mathcal{C}}
\def\O{{\mathcal{O}}}

\def\cc{{\bf c}}
\def\ee{{\bf e}}
\def\x{{\bf x}}
\def\y{{\bf y}}

\def\z{{\bf z}}

\def\uu{{\bf u}}
\def\vv{{\bf v}}
\def\0{{\bf 0}}

\title{\bf{Tensor Norm, Cubic Power and Gelfand Limit}}%\thanks{This research was supported by the Hong Kong
\author{ \hspace{1mm} Liqun Qi\thanks{Department of Applied
    Mathematics, The Hong Kong Polytechnic University, Hung Hom,
    Kowloon, Hong Kong; ({\tt liqun.qi@polyu.edu.hk}). This author's work was supported by the Hong Kong
    Research Grant Council (Grant No.  PolyU 15300715, 15301716 and 15300717). },
    \ \
    Shenglong Hu\thanks{Department of Mathematics, School of Science, Hangzhou Dianzi University, Hangzhou 310018 China; ({\tt shenglonghu@hdu.edu.cn}). This author's work was supported by NSFC (Grant No.  11771328).},
 \ \
    Xinzhen Zhang\thanks{School of Mathematics, Tianjin University, Tianjin 300354 China; ({\tt xzzhang@tju.edu.cn}). This author's work was supported by NSFC (Grant No.  11871369). }
 \ and \
    Yannan Chen\thanks{School of Mathematical Sciences, South China    Normal University, Guangzhou, China; ({\tt ynchen@scnu.edu.cn}).  This author was supported by the National Natural Science Foundation of China (11771405).}}

\begin{document}
\date{\today}
\maketitle

\begin{abstract}
We establish two inequalities for the nuclear norm and the spectral norm of tensor products.  The first inequality indicates that the nuclear norm of the square matrix is a matrix norm.  We extend the concept of matrix norm to tensor norm.  We show that the $1$-norm, the Frobenius norm and the nuclear norm of tensors are tensor norms, but the infinity norm and the spectral norm of tensors are not tensor norms.   We introduce the cubic power for a general third order tensor, and show that a Gelfand formula holds for a general third order tensor.
In that formula, for any norm, a common spectral radius-like limit exists for that third order tensor.   We call such a limit the Gelfand limit. The Gelfand limit is zero if the third order tensor is nilpotent, and is one or zero if the third order tensor is idempotent.   The Gelfand limit is not greater than any tensor norm of that third order tensor, and the cubic power of that third order tensor tends to zero as the power increases to infinity if and only if the Gelfand limit is less than one.   The cubic power and the Gelfand limit can be extended to any higher odd order tensors.

%Finally, we apply our results to biquadratic tensors.


%15A69:LINEAR AND MULTILINEAR ALGEBRA; MATRIX THEORY- Multilinear algebra, tensor products
%53A45:Classical differential geometry-Vector and tensor analysis
%47A05: Operator Theory-General (adjoints, conjugates, products, inverses, domains, ranges,
%etc.)
%53C35:Global differential geometry- Symmetric spaces
\end{abstract}

%\newpage

\section{Introduction}

The tensor completion and recovery problem aims to fill the missing or unobserved entries of partially observed tensors.  It has received wide attention and achievements in areas like data mining, computer vision, signal processing, and neuroscience \cite{SGCH19, YZ16}.   Former approaches often proceed by unfolding tensors to matrices and then apply for matrix completion.
Yuan and Zhang \cite{YZ16} showed that such matricization fails to exploit the structure of tensors and may lead to sub-optimality.  They proposed to minimize a tensor nuclear norm directly and proved that such an approach improves the sample size requirement.    This leads research enthusiasm on the tensor nuclear norm and its dual norm, i.e., the tensor spectral norm \cite{Hu15, JYZ17, KLW18, Li16, Ni17}, though this is in fact a NP-hard problem \cite{FL17}.

In matrix analysis \cite{HJ12}, the matrix norm is a concept different from the vector norm.  We may regard a matrix space as a vector space.   A norm on that matrix space is called a vector norm. If in additional, that norm satisfies the axiom that the norm of the product of two arbitrary matrices is always not greater than the product of the norms of these two matrices, then that norm is called a matrix norm.   It was shown that the $1$-norm, the $2$-norm and all the induced norms of matrices are matrix norms, but the infinity norm of matrices is not a matrix norm \cite{HJ12}.

 The matrix norm of a square matrix $A$ is closely linked with the spectral radius of $A$, i.e., $\rho(A)$ , by the well-known Gelfand formula (1941):
$$\rho(A)  = \lim_{k\to\infty} |||A^k|||^{1 \over k},$$
for any matrix norm $|||\cdot |||$.

However, there is no ``tensor norm'' concept or a Gelfand formula for higher order tensors until now.   In this paper, we explore this unknown territory.

In the next section, we show that the nuclear norm of the tensor product of two tensors is not greater than the product of the nuclear norms of these two tensors.

However, in general, the spectral norm of the tensor product of two tensors may be greater than the product of the spectral norms of these two tensors.   In Section 3, we give a counterexample to illustrate this.   As a substitute, in that section, we show that the spectral norm of the tensor product of two tensors is not greater than the product of the spectral norm of one tensor, and the nuclear norm of another tensor.

 By the result in Section 2, we conclude in Section 4 that the nuclear norm of the square matrix is also a matrix norm.   Viewing the significance of the nuclear norm of matrices in the matrix completion problem \cite{CR08}, and the importance of the matrix norm in matrix analysis \cite{HJ12}, this result may be useful in the related research.

In Section 5, we extend the concept of matrix norm to tensor norm.  A real function defined for all real tensors is called a tensor norm if it is a norm for any tensor space with fixed dimensions, and the norm of the tensor product of two tensors is always not greater than the product of the norms of these two tensors.  We show that the $1$-norm, the Frobenius norm and the nuclear norm of tensors are tensor norms but the infinity norm and the spectral norm of tensors are not tensor norms.

In Section 6, we introduce the cubic power for a general third order tensor.   The cubic power preserves nonnegativity, symmetry and the diagonal property.   We show that the cubic power of a general third order tensor tends to zero as the power increases to infinity, if there is a tensor norm such that the tensor norm of that third order tensor is less than one.

We show in Section 7 that a Gelfand formula holds for a general third order tensor.
In that formula, for any norm, a common spectral radius-like limit exists for that third order tensor.   We call such a limit the Gelfand limit of that third order tensor. The Gelfand limit is zero if the  third order tensor is nilpotent, and is one or zero if the third order tensor is idempotent. We show that the Gelfand limit is less than or equal to any tensor norm of that third order tensor, and the cubic power of that third order tensor tends to zero as the power increases to infinity if and only if the Gelfand limit is less than one.

The cubic power and the Gelfand limit can be extended to any higher odd order tensors.

%Finally, we apply our results to biquadratic tensors.

Some final remarks are made in Section 8.

In Sections 2, 3, 6 and 7, there is a subsection for each section.  These subsections describe some sub-developments of such sections.  The reader may skip them in the first reading.

In this paper, unless otherwise stated, all the discussions will be carried out in the filed of real numbers.   We use small letters $\lambda, x_i, u_i$, etc., to denote scalars, small bold letters $\x, \uu, \vv$, etc.,  to denote vectors, capital letters $A, B, C$, etc., to denote matrices, and calligraphic letters $\A, \B, \C$, etc., to denote tensors, with $\O$ as the zero tensor with adequate order and dimensions.

\section{Nuclear Norm of Tensor Product}

Let $\mathbb{N}$ be the set of positive integers, and $\bar {\mathbb{N}}$ be the set of nonnegative integers.   For $k \in \mathbb{N}$, we use $[k]$ to denote the set $\{ 1, \cdots, k\}$.   For a vector $\uu = (u_1, \cdots, u_n)^\top$, we use $\|\uu\|_2$ to denote its 2-norm.   Thus,
$$\| \uu\|_2 := \sqrt{u_1^2 + \cdots + u_n^2}.$$

Suppose that a $k$th order tensor $\A = (a_{i_1\cdots i_k}) \in \Re^{n_1 \times \cdots \times n_k}$, where $k \in \mathbb{N}$ is called the order of $\A$, and $n_i \in \mathbb{N}$ for $i \in [k]$ are called the dimensions of $\A$.    We use $\circ$ to denote tensor outer product.  Then for nonzero $\uu^{(i)} \in \Re^{n_i}, i \in [k]$,
$$\uu^{(1)} \circ \cdots \circ \uu^{(k)}$$
is a rank-one $k$th order tensor.
The nuclear norm of $\A$ is defined \cite{FL17, Hu15, Li16} as
\begin{equation} \label{e1}
\| \A \|_* := \inf \left\{ \sum_{j=1}^r |\lambda_j| : \A = \sum_{j=1}^r \lambda_j \uu^{(1, j)} \circ \cdots \circ \uu^{(k, j)}, \left\|\uu^{(i, j)}\right\|_2 = 1, \uu^{(i, j)} \in \Re^{n_i}, i \in [k], j \in [r], r \in \mathbb{N}
\right\}.
\end{equation}

We have the following theorem.
\begin{theorem} \label{t1}
Suppose that $\A = (a_{i_1\cdots i_{k+p}}) \in \Re^{n_1 \times \cdots \times n_{k+p}}$, $\B = (b_{i_{k+1}\cdots i_{k+p+q}}) \in \Re^{n_{k+1} \times \cdots \times n_{k+p+q}}$, and a tensor product of $\A$ and $\B$ is defined as $\C = (c_{i_1\cdots i_ki_{k+p+1}\cdots i_{k+p+q}}) \in \Re^{n_1 \times \cdots \times n_k \times n_{k+p+1} \cdots \times n_{k+p+q}}$ by
$$c_{i_1\cdots i_ki_{k+p+1}\cdots i_{k+p+q}} = \sum_{i_{k+1}=1}^{n_{k+1}} \cdots \sum_{i_{k+p}=1}^{n_{k+p}}a_{i_1\cdots i_{k+p}}b_{i_{k+1}\cdots i_{k+p+q}},$$
for $i_l \in [n_l],$ $l = 1, \cdots, k, k+p+1, \cdots, k+p+q$, with $k, p\in \mathbb{N}$, $q \in \bar {\mathbb{N}}$.
Then
\begin{equation} \label{e2}
\| \C \|_* \le \|\A \|_* \|\B \|_*.
\end{equation}
\end{theorem}
{\bf Proof} Let $\epsilon > 0$.   Then we have
$$\A = \sum_{j=1}^{r_1} \lambda_j \uu^{(1, j)} \circ \cdots \circ \uu^{(k+p, j)},$$
where $\|\uu^{(i, j)}\| = 1$ for $i = 1, \cdots, k+p$, and $j \in [r_1], r_1 \in \mathbb{N}$, and
$$\B = \sum_{l=1}^{r_2} \mu_l \vv^{(k+1, l)} \circ \cdots \circ \vv^{(k+p+q, l)},$$
where $\|\vv^{(i, l)}\| = 1$ for $i = k+1, \cdots, k+p+q$, and $l \in [r_2], r_2 \in \mathbb{N}$,
such that
$$\sum_{j=1}^{r_1} |\lambda_j| \le \| \A \|_* + \epsilon$$
and
$$\sum_{l=1}^{r_2} |\mu_l| \le \| \B \|_* + \epsilon.$$
This implies that
$$\C = \sum_{j=1}^{r_1} \sum_{l=1}^{r_2} \lambda_j\mu_l \prod_{i=1}^p \langle \uu^{(k+i, j)}, \vv^{(k+i, l)} \rangle \uu^{(1, j)} \circ \cdots \circ \uu^{(k, j)}\circ \vv^{(k+p+1, l)} \circ \cdots \circ \vv^{(k+p+r, l)}.$$
Hence,
$$\begin{aligned} \| \C \|_* & \le \sum_{j=1}^{r_1} \sum_{l=1}^{r_2} |\lambda_j\mu_l \prod_{i=1}^p \langle \uu^{(k+i, j)}, \vv^{(k+i, l)} \rangle|\\
& \le \sum_{j=1}^{r_1} \sum_{l=1}^{r_2} |\lambda_j\mu_l| \prod_{i=1}^p \| \uu^{(k+i, j)}\|_2 \|\vv^{(k+i, l)}\|_2\\
& = \sum_{j=1}^{r_1} \sum_{l=1}^{r_2} |\lambda_j\mu_l|\\
& = \left(\sum_{j=1}^{r_1} |\lambda_j|\right)\left(\sum_{l=1}^{r_2} |\mu_l|\right)\\
& \le \left( \|\A \|_* + \epsilon \right)\left( \|\B \|_* + \epsilon \right),
\end{aligned}$$
where the first inequality is by the definition (\ref{e1}), and the second inequality is by the Cauchy inequality.
Letting $\epsilon \to 0$, we have (\ref{e2}).
\qed

%The remaining part of this section is an application of this theorem, and is not related with the main part of this paper.  A reader who is interested in tensor norm, cubic power and Gelfand limit may skip the remaining part of this section.

\subsection{An Application: Lower Bounds for the Nuclear Norm of a Tensor}

In this subsection, we present a lower bound for the nuclear norm of an arbitrary even order tensor $\A = (a_{i_1\cdots i_k}) \in \Re^{n_1 \times \cdots \times n_k}$ for $k \in \mathbb{N}$ and $k \ge 3$.

We first assume that $k=3$.   Then $\A$ is a third order tensor.   It is not easy to compute its nuclear norm.

\begin{proposition} \label{p3.1}
 Suppose that $\A = (a_{i_1i_2}) \in \Re^{n_1 \times n_2 \times n_3}$, where $n_1, n_2, n_3 \in \mathbb{N}$.  Let $B = (b_{i_2i_3}) \in \Re^{n_2 \times n_3}$, and $\cc  = (c_1, \cdots, c_{n_1})^\top \in \Re^{n_1}$ be defined by
 \begin{equation} \label{e3.3}
 c_i = \sum_{i_2 = 1}^{n_2} \sum_{i_3 = 1}^{n_3} a_{ii_2i_3}b_{i_2i_3},
 \end{equation}
 for $i \in [n_1]$.  Then,
 \begin{equation} \label{e3.4}
 \|\A \|_* \ge  \max \left\{ \| \cc \|_* : \cc \ {\rm is \ calculated \ by \ (\ref{e3.3})}, B \in \Re^{n_2 \times n_3}, \| B \|_* = 1 \right\}.
 \end{equation}
\end{proposition}
{\bf Proof}  Applying Theorem \ref{t1} with $\B = B$ and $\C = \cc$, $k = 1$, $p = 2$ and $q = 0$, we have the conclusion.
\qed

Here, $B$ is a matrix, and $\cc$ is a vector.   Their nuclear norms are not difficult to be calculated.  The tensor $\A$ in the following example is originally from \cite{FL17}.
\begin{example}
Let $k = 3$, $n_1 = n_2 = n_3 = 2$, and $\A$ be a third order symmetric tensor defined by
$$\A = {1 \over 2}\left( \ee^{(1)} \circ \ee^{(1)} \circ \ee^{(2)} + \ee^{(1)} \circ \ee^{(2)} \circ \ee^{(1)} + \ee^{(2)} \circ \ee^{(1)} \circ \ee^{(1)} - \ee^{(2)} \circ \ee^{(2)} \circ \ee^{(2)}\right),$$
where $\ee^{(1)} = (1, 0)^\top$ and $\ee^{(2)} = (0, 1)^\top$.  Let $J$ be the all one matrix in $\Re^{2 \times 2}$, and
$B = {J \over \|J\|_*}$.  We may calculate $\cc$ by (\ref{e3.4}).  Then by (\ref{e3.4}), we have
$$\|\A\|_* \ge \|\cc\|_* \equiv 0.6455.$$
Actually, $\|\A \|_* = 2$.   This verifies (\ref{e3.4}) somehow.
\end{example}

We then consider the case that $k = 4$.  Then $\A$ is a fourth order tensor.   It is also not easy to compute its nuclear norm.

\begin{proposition} \label{p3.3}
 Suppose that $\A = (a_{i_1i_2i_3i_4}) \in \Re^{n_1 \times n_2 \times n_3 \times n_4}$, where $n_1, n_2, n_3, n_4 \in \mathbb{N}$.  Let $B = (b_{i_3i_4}) \in \Re^{n_3 \times n_4}$, and $C = (c_{i_1i_2}) \in \Re^{n_1 \times n_2}$ be defined by
 \begin{equation} \label{e3.5}
 c_{i_1i_2} = \sum_{i_3 = 1}^{n_3} \sum_{i_4 = 1}^{n_4} a_{i_1i_2i_3i_4}b_{i_3i_4},
 \end{equation}
 for $i_1 \in [n_1], i_2 \in [n_2]$.  Then,
 \begin{equation} \label{e3.6}
 \|\A \|_* \ge  \max \left\{ \| C \|_* : C \ {\rm is \ calculated \ by \ (\ref{e3.5})}, B \in \Re^{n_3 \times n_4}, \| B \|_* = 1 \right\}.
 \end{equation}
\end{proposition}
{\bf Proof}  Applying Theorem \ref{t1} with $\B = B$ and $\C = C$, $k = p = 2$ and $q = 0$, we have the conclusion.
\qed

Here, $B$ and $C$ are matrices.   Their nuclear norms are not difficult to be calculated.  The following example is from \cite{Ni17}.

\begin{example}
Let $k=4$, $n_1 = n_2 = n_3 = n_4 = 3$, and $\A$ be a fourth order symmetric tensor defined by
$$\A = \ee^{\otimes 4} - \left(\ee^{(1)}\right)^{\otimes 4} - \left(\ee^{(2)}\right)^{\otimes 4} - \left(\ee^{(3)}\right)^{\otimes 4},$$
where $\ee = (1, 1, 1)^\top$, $\ee^{(1)} = (1, 0, 0)^\top$, $\ee^{(2)} = (0, 1, 0)^\top$ and $\ee^{(3)} = (0, 0, 1)^\top$.
Let $J$ be the all one matrix in $\Re^{3 \times 3}$, and
$B = {J \over \|J\|_*}$.   We may calculate $C$ by (\ref{e3.5}).  Then by (\ref{e3.6}), we have
$$\|\A\|_* \ge \|C\|_* \equiv 10.3757.$$
Actually, $\|\A \|_* = 12$.   This verifies (\ref{e3.6}) somehow.
\end{example}

We now extend Propositions \ref{p3.1} and \ref{p3.3} to the case that $k \ge 5$.   We have the following two propositions.

\begin{proposition} \label{p3.5}
 Suppose that $\A = (a_{i_1\cdots i_k}) \in \Re^{n_1 \times \cdots \times n_k}$, where $k = 2l+1$, $n_1, \cdots, n_k, l \in \mathbb{N}, l \ge 2$.  Let $B^{(j)} = \left(b^{(j)}_{i_{2j}i_{2j+1}}\right) \in \Re^{n_{2j} \times n_{2j+1}}$ for $j = 1, \cdots, l$, and $\cc  = (c_1, \cdots, c_{n_1})^\top \in \Re^{n_1}$ be defined by
 \begin{equation} \label{e3.7}
 c_i = \sum_{i_2 = 1}^{n_2} \cdots \sum_{i_k = 1}^{n_k} a_{ii_2\cdots i_k}b^{(1)}_{i_2i_3}\cdots b^{(l)}_{i_{2l}i_{2l+1}},
 \end{equation}
 for $i \in [n_1]$.  Then,
 \begin{equation} \label{e3.8}
 \|\A \|_* \ge  \max \left\{ \| \cc \|_* : \cc \ {\rm is \ calculated \ by \ (\ref{e3.7})}, \B^{(j)} \in \Re^{n_{2j} \times n_{2j+1}}, \| B^{(j)} \|_* = 1, j \in [l] \right\}.
 \end{equation}
\end{proposition}
{\bf Proof}  Applying Theorem \ref{t1} repetitively, we have the conclusion.
\qed

Here, $B^{(j)}$ for $j \in [l]$ are matrices, and $\cc$ is a vector.   Their nuclear norms are not difficult to be calculated.

\begin{proposition} \label{p3.6}
 Suppose that $\A = (a_{i_1\cdots i_k}) \in \Re^{n_1 \times \cdots \times n_k}$, where $k = 2l+2$, $n_1, \cdots, n_k, l \in \mathbb{N}, l \ge 2$.  Let $B^{(j)} = (b^{(j)}_{i_{2j+1}i_{2j+2}}) \in \Re^{n_{2j+1} \times n_{2j+2}}$ for $j = 1, \cdots, l$, and $C  = (c_{i_1i_2})^\top \in \Re^{n_1 \times n_2}$ be defined by
 \begin{equation} \label{e3.7}
 c_{i_1i_2} = \sum_{i_3 = 1}^{n_3} \cdots \sum_{i_k = 1}^{n_k} a_{i_1\cdots i_k}b^{(1)}_{i_3i_4}\cdots b^{(l)}_{i_{2l+1}i_{2l+2}},
 \end{equation}
 for $i_1 \in [n_1]$ and $i_2 \in [n_2]$.  Then,
 \begin{equation} \label{e3.8}
 \|\A \|_* \ge  \max \left\{ \| C \|_* : C \ {\rm is \ calculated \ by \ (\ref{e3.7})}, \B^{(j)} \in \Re^{n_{2j+1} \times n_{2j+2}}, \| B^{(j)} \|_* = 1, j \in [l] \right\}.
 \end{equation}
\end{proposition}
{\bf Proof}  Applying Theorem \ref{t1} repetitively, we have the conclusion.
\qed

Here, $B^{(j)}$ for $j \in [l]$, and $C$ are matrices.   Their nuclear norms are not difficult to be calculated.

\section{Spectral Norm of Tensor Product}

For $\A = (a_{i_1\cdots i_k}), \B = (b_{i_1\cdots i_k}) \in \Re^{n_1 \times \cdots \times n_k}$, their inner product is defined as
$$\langle \A, \B \rangle := \sum_{i_1=1}^{n_1} \cdots  \sum_{i_k=1}^{n_k} a_{i_1\cdots i_k}b_{i_1\cdots i_k}.$$
Then the spectral norm of $\A$ is defined \cite{FL17, Hu15, Li16, YZ16} as
\begin{equation} \label{e3}
\| \A \|_S := \max \left\{ \langle \A, \uu^{(1)} \circ \cdots \circ \uu^{(k)} \rangle, \uu^{(i)} \in \Re^{n_i}, \left\|\uu^{(i)}\right\|_2 = 1, \ {\rm for}\ i \in [k] \right\}.
\end{equation}
It is known \cite{YZ16} that we always have
$$\|\A \|_S \le \|\A \|_*.$$
Note that in general, we do not have
$$\|\C \|_S \le \|\A \|_S \|\B \|_S,$$
if $\C$ is a tensor product of $\A$ and $\B$, as in Theorem \ref{t1}.   See the following example.

\begin{example}  \label{ex3.1}
Let $k = p = q = 2$, and $n_i = 2$ for $i =1, \cdots, 6$.   Let $\A = (a_{i_1i_2i_3i_4}) = \B = (b_{i_3i_4i_5i_6})$ be defined by
$$\begin{aligned}
a_{1111} = 2, \ \ & a_{1211} = 3, & a_{2111} = -6, \ \ & a_{2211} = 3,\\
a_{1121} = -6, \ \ & a_{1221} = 3, & a_{2121} = 4, \ \ & a_{2221} = 3,\\
a_{1112} = 3, \ \ & a_{1212} = 9, & a_{2112} = 3, \ \ & a_{2212} = -3,\\
a_{1122} = 3, \ \ & a_{1222} = -3, & a_{2122} = 3, \ \ & a_{2222} = 15.
\end{aligned}$$
Then we have
$$\begin{aligned}
c_{1111} = 58, \ \ & c_{1211} = 6, & c_{2111} = -18, \ \ & c_{2211} = 24,\\
c_{1121} = -18, \ \ & c_{1221} = 12, & c_{2121} = 70, \ \ & c_{2221} = 30,\\
c_{1112} = 6, \ \ & c_{1212} = 108, & c_{2112} = 12, \ \ & c_{2212} = -54,\\
c_{1122} = 24, \ \ & c_{1222} = -54, & c_{2122} = 30, \ \ & c_{2222} = 1252.
\end{aligned}$$
By calculation, we have $\| \A \|_S = \|\B \|_S = 16.3609$, and $\| \C \|_S = 271.5503$.  Then $\| \A\|_S \|\B \|_S = 268.6781$, which is slightly less than $\| \C \|_S$.
\end{example}

However, we may establish the following theorem.

\begin{theorem} \label{t2}
Suppose that $\A = (a_{i_1\cdots i_{k+p}}) \in \Re^{n_1 \times \cdots \times n_{k+p}}$, $\B = (b_{i_{k+1}\cdots i_{k+p+q}}) \in \Re^{n_{k+1} \times \cdots \times n_{k+p+q}}$, and a tensor product of $\A$ and $\B$ is defined as $\C = (c_{i_1\cdots i_ki_{k+p+1}\cdots i_{k+p+q}}) \in \Re^{n_1 \times \cdots \times n_k \times n_{k+p+1} \cdots \times n_{k+p+q}}$ by
$$c_{i_1\cdots i_ki_{k+p+1}\cdots i_{k+p+q}} = \sum_{i_{k+1}=1}^{n_{k+1}} \cdots \sum_{i_{k+p}=1}^{n_{k+p}}a_{i_1\cdots i_{k+p}}b_{i_{k+1}\cdots i_{k+p+q}},$$
for $i_l \in [n_l],$ $l = 1, \cdots, k, k+p+1, \cdots, k+p+q$, with $k, p \in \mathbb{N}$, $q \in \bar {\mathbb{N}}$.
Then
\begin{equation} \label{e4}
\| \C \|_S \le \|\A \|_S \|\B \|_*.
\end{equation}
\end{theorem}
{\bf Proof} Let $\epsilon > 0$.   Then we have
$$\B = \sum_{j=1}^r \mu_j \vv^{(k+1, j)} \circ \cdots \circ \vv^{(k+p+q, j)},$$
where $\|\vv^{(i, j)}\| = 1$ for $i = k+1, \cdots, k+p+q$, and $j \in [r], r \in \mathbb{N}$,
such that
$$\sum_{j=1}^r |\mu_j| \le \| \B \|_* + \epsilon.$$
Then
$$b_{i_{k+1}\cdots i_{k+p+q}} = \sum_{j=1}^r \mu_jv^{(k+1, j)}_{i_{k+1}} \cdots v^{(k+p+q, j)}_{i_{k+p+q}},$$
for $i = k+1, \cdots, k+p+q$, and $j \in [r]$.

We have
$$c_{i_1\cdots i_ki_{k+p+1}\cdots i_{k+p+r}} = \sum_{j=1}^r \sum_{i_{k+1}=1}^{n_{k+1}} \cdots \sum_{i_{k+p}=1}^{n_{k+p}} \mu_ja_{i_1\cdots i_{k+p}} v^{(k+1, j)}_{i_{k+1}} \cdots v^{(k+p+q, j)}_{i_{k+p+q}},$$
for $i_l \in [n_l],$ $l = 1, \cdots, k, k+p+1, \cdots, k+p+q$, and $j \in [r]$.
This implies that
$$\begin{aligned}
& \|\C\|_S \\
= & \max \left\{ \left\langle \C, \uu^{(1)} \circ \cdots \circ \uu^{(k)} \circ \uu^{(k+p+1)} \circ \cdots \circ \uu^{(k+p+q)} \right\rangle : \| \uu^{(i)}\|_2 = 1 \right\}\\
= & \max \left\{ \sum_{j=1}^r \mu_j \sum_{i_1=1}^{n_1} \cdots \sum_{i_{k+p+q=1}}^{n_{k+p+q}}a_{i_1\cdots i_{k+p}} u_{i_1}^{(1)} \cdots u_{i_k}^{(k)}v^{(k+1, j)}_{i_{k+1}} \cdots v^{(k+p+q, j)}_{i_{k+p+q}}u_{i_{l+p+1}}^{(k+p+1)} \cdots u_{i_{k+p+q}}^{(k+p+q)} \right\}\\
 \le & \max \left\{ \sum_{j=1}^r \left|\mu_j\right| \left|\left\langle \A, \uu^{(1)} \circ \cdots \circ \uu^{(k)} \circ \vv^{(k+1, j)} \circ \cdots \circ \vv^{(k+p, j)} \right\rangle\right| \prod_{l=k+p+1}^{k+p+q} \left|\left\langle \vv^{(l, j)}, \uu^{(l)} \right\rangle\right| \right\}\\
 \le & \max \left\{ \sum_{j=1}^r \left|\mu_j\right| \|\A \|_S \right\}\\
 \le & \|\A \|_S \left(\|\B \|_* + \epsilon \right),
\end{aligned}$$
where the second inequality is by the definition (\ref{e3}) and the Cauchy inequality.
Letting $\epsilon \to 0$, we have (\ref{e4}).
\qed

%The remaining part of this section contains some applications of this theorem, and is not related with the main part of this paper.  A reader who is interested in tensor norm, cubic power and Gelfand limit may skip the remaining part of this section.

\subsection{An Alternative Formula for the Spectral Norm of a Tensor}

We present an alternative formula for the spectral norm of a tensor in this subsection.  This formula does not reduce the complexity of the problem, as this is impossible, but gives an alternative approach to handle the spectral norm.

\begin{proposition}
Suppose that $\A = (a_{i_1\cdots i_k}) \in \Re^{n_1 \times \cdots \times n_k}$, where $k \in \mathbb{N}$, $k \ge 3$.  For $B = (b_{i_{k-1}i_k}) \in \Re^{n_{k-1} \times n_k}$, define $\C = (c_{i_1\cdots i_{k-2}}) \in \Re^{n_1 \times \cdots \times n_{k-2}}$ by
\begin{equation} \label{ea.1}
c_{i_1\cdots i_{k-2}} = \sum_{i_{k-1}=1}^{n_{k-1}} \sum_{i_k=1}^{n_k}a_{i_1\cdots i_k}b_{i_{k-1}i_k},
\end{equation}
for $i_j \in [n_j], j = 1, \cdots, k-2$.
Then
\begin{equation} \label{ea.2}
\|\A \|_S = \max \left\{ \| \C \|_S : \C \ {\rm is \ calculated \ by \ (\ref{ea.1})}, B \in \Re^{n_{k-1} \times n_k}, \|B\|_* = 1 \right\}.
\end{equation}
\end{proposition}
{\bf Proof} By Theorem \ref{t2}, we have
$$\|\A \|_S \ge \max \left\{ \| \C \|_S : \C \ {\rm is \ calculated \ by \ (\ref{ea.1})}, B \in \Re^{n_{k-1} \times n_k}, \|B\|_* = 1 \right\}.$$
On the other hand,
$$\begin{aligned} & \| \A \|_S \\
& = \max \left\{ \langle \A, \uu^{(1)} \circ \cdots \circ \uu^{(k)} \rangle, \uu^{(i)} \in \Re^{n_i}, \left\|\uu^{(i)}\right\|_2 = 1, \ {\rm for}\ i \in [k] \right\} \\
& \le \max \left\{ \langle \A, \uu^{(1)} \circ \cdots \circ \uu^{(k-2)} \circ B \rangle, \uu^{(i)} \in \Re^{n_i}, \left\|\uu^{(i)}\right\|_2 = 1, \ {\rm for}\ i \in [k-2], B \in \Re^{n_{k-1} \times n_k}, \|B\|_* = 1 \right\} \\
& = \max \left\{ \| \C \|_S : \C \ {\rm is \ calculated \ by \ (\ref{ea.1})}, B \in \Re^{n_{k-1} \times n_k}, \|B\|_* = 1 \right\}.
\end{aligned}$$
This proves (\ref{ea.2}).
\qed

\subsection{Lower Bounds for the Product of the Nuclear Norm and Spectral Norm of a Tensor}

We now present some lower bounds for the product of the nuclear norm and the spectral norm of $\A \in \Re^{n_1 \times \cdots \times n_k}$.   Denote the spectral radius of a matrix $A \in \Re^{n \times n}$ by $\rho(A)$.  We introduce contraction matrices of $\A$.

\begin{definition}
Let $\A = (a_{i_1\cdots i_k}) \in \Re^{n_1 \times \cdots \times n_k}$, where $n_1, \cdots, n_k, k \in \mathbb{N}$. Assume that $k \ge 2$.   Let $j \in [k]$.  Define a symmetric matrix $A^{(j)} = \left(a^{(j)}_{rs}\right) \in \Re^{n_j \times n_j}$ by
$$a^{(j)}_{rs} = \sum_{i_1=1}^{n_1} \cdots \sum_{i_{j-1}=1}^{n_{j-1}}\sum_{i_{j+1}=1}^{n_{j+1}}\cdots \sum_{i_k=1}^{n_k} a_{i_1\cdots i_{k-1}ri_{k+1}\cdots i_k}a_{i_1\cdots i_{k-1}si_{k+1}\cdots i_k},$$
for $r, s \in [n_j]$.  We call $A^{(j)}$ the $j$th contraction matrix of $\A$.
\end{definition}

\begin{proposition} \label{p5.1}
 Suppose that $\A = (a_{i_1\cdots i_k}) \in \Re^{n_1 \times \cdots \times n_k}$, where $n_1, \cdots, n_k, k \in \mathbb{N}$. Assume that $k \ge 2$.  Let $j \in [k]$.   Then,
 \begin{equation} \label{e5.6}
 \rho(A^{(j)}) \equiv \|A^{(j)}\|_S \le \|\A \|_*\|\A\|_S.
 \end{equation}
\end{proposition}
{\bf Proof}  Apply Theorem \ref{t2} with $\B = \A$ and $\C = A^{(j)}$.  Note that $A^{(j)} \in \Re^{n_j \times n_j}$ is symmetric.  Hence, its spectral norm is its spectral radius.   We then have the conclusion.
\qed

As $A^{(j)}$ is a symmetric matrix, its spectral norm, i.e., its spectral radius is the largest absolute value of its eigenvalues, which is not difficult to be calculated.   For $j = 1, \cdots , k$, this theorem gives $k$ lower bounds for the product of the nuclear norm and the spectral norm of $\A$.   The tensor $\A$ in the following example is the same as the tensor $\A$ in Example 2.3.  It is originally from \cite{FL17}.
\begin{example}
Let $k = 3$, $n_1 = n_2 = n_3 = 2$, and $\A$ be a third order symmetric tensor defined by
$$\A = {1 \over 2}\left( \ee^{(1)} \circ \ee^{(1)} \circ \ee^{(2)} + \ee^{(1)} \circ \ee^{(2)} \circ \ee^{(1)} + \ee^{(2)} \circ \ee^{(1)} \circ \ee^{(1)} - \ee^{(2)} \circ \ee^{(2)} \circ \ee^{(2)}\right),$$
where $\ee^{(1)} = (1, 0)^\top$ and $\ee^{(2)} = (0, 1)^\top$.
Then
$$A^{(1)} = A^{(2)} = A^{(3)} = 0.5I_2,$$
where $I_2$ is the identity matrix in $\Re^{2 \times 2}$.   Then we have $\|\A \|_S = 0.5$, $\|\A\|_* = 2$ and $\rho(A^{(1)}) = \rho(A^{(2)}) = \rho(A^{(3)}) = 0.5$.  This verifies (\ref{e5.6}).
\end{example}

\section{Matrix Norm}

In matrix analysis \cite{HJ12}, the matrix norm is a concept different from the vector norm.   Consider matrices in $\Re^{n \times n}$.   Let
$||| \cdot ||| : \Re^{n \times n} \to \Re_+$.   If it is not only a vector norm in $\Re^{n \times n}$, but it also satisfies the following additional axiom: for any $A, B \in \Re^{n \times n}$,
\begin{equation} \label{e7.14}
||| A B ||| \le ||| A ||| \cdot ||| B |||,
\end{equation}
then $||| \cdot |||$ is a matrix norm in $\Re^{n \times n}$.   Otherwise, it is only a vector norm.   In particular, $1$-norm, $2$-norm and any induced norm are matrix norms, but $\infty$-norm is only a vector norm, not a matrix norm.   Matrix norms play an important role in matrix analysis.  See \cite{HJ12} for more details on matrix norms.

By Theorem \ref{t1}, we have the following proposition.

\begin{proposition} \label{p1}
For $\Re^{n \times n}$, the nuclear norm is a matrix norm.
\end{proposition}

As the nuclear norm plays a significant role in the matrix completion problem \cite{CR08}, this conclusion for matrix nuclear norm should be useful in the related research.

By Theorems \ref{t1} and \ref{t2}, we have the following further results.

\begin{proposition} \label{p2}
Suppose that $A \in \Re^{n \times n}$ is invertible.   Then we have
\begin{equation} \label{e5}
\|A\|_*\|A^{-1}\|_*  \ge n,
\end{equation}
and
\begin{equation} \label{e6}
\|A\|_*\|A^{-1}\|_S  \ge 1.
\end{equation}
\end{proposition}
{\bf Proof} Apply (\ref{e2}) and (\ref{e4}) to $A$ and $\A^{-1}$.   Note that $\| I_n \|_S = 1$ and $\| I_n \|_* = n$, where $I_n$ is the identity matrix in $\Re^{n \times n}$.   The conclusions hold.
\qed

\begin{proposition} \label{p3}
Suppose that $A \in \Re^{n \times n}$.   If $\|A\|_* < 1$, then
$$ \lim_{k \to \infty} A^k = 0.$$
\end{proposition}
{\bf Proof} If $\|A\|_* < 1$, then
$$\|A^k\|_* \le \|A\|_*^k \to 0.$$
The result follows. \qed

\section{Tensor Norm}

We are now ready to extend the concept of matrix norms to tensor norms.

\begin{definition}
Suppose that $||| \cdot |||$ is a function defined for all real tensors, and in any real tensor space $\Re^{n_1 \times \cdots \times n_k}$ of fixed dimensions with $n_1, \cdots, n_k, k \in \mathbb{N}$, it is a vector norm.  If furthermore for any two real tensors $\A$ and $\B$ such that $\A$ and $\B$ have an outer tensor product $\C$, we have
\begin{equation} \label{e7}
|||\C ||| \le ||| \A ||| \cdot ||| \B |||,
\end{equation}
then $||| \cdot |||$ is called a tensor norm.
\end{definition}

Clearly, a tensor norm must be a matrix norm if it is restricted to $\Re^{n \times n}$.  We have the following theorem.

\begin{theorem} \label{t3}
The nuclear norm, the 1-norm, the Frobenius norm are tensor norms, but the infinity norm and the spectral norm are not tensor norms.
\end{theorem}
{\bf Proof}  By Theorem \ref{t1}, the nuclear norm is a tensor norm.   Since the infinity norm is not a matrix norm, it is also not a tensor norm.  By the counter example in Example 3.1, the spectral norm is not a tensor norm.  What we need to check are the 1-norm and the Frobenius norm.

Suppose that $\A = (a_{i_1\cdots i_{k+p}}) \in \Re^{n_1 \times \cdots \times n_{k+p}}$, $\B = (b_{i_{k+1}\cdots i_{k+p+q}}) \in \Re^{n_{k+1} \times \cdots \times n_{k+p+q}}$, and a tensor product of $\A$ and $\B$ is defined as $\C = (c_{i_1\cdots i_ki_{k+p+1}\cdots i_{k+p+q}}) \in \Re^{n_1 \times \cdots \times n_k \times n_{k+p+1} \cdots \times n_{k+p+q}}$ by
$$c_{i_1\cdots i_ki_{k+p+1}\cdots i_{k+p+q}} = \sum_{i_{k+1}=1}^{n_{k+1}} \cdots \sum_{i_{k+p}=1}^{n_{k+p}}a_{i_1\cdots i_{k+p}}b_{i_{k+1}\cdots i_{k+p+q}},$$
for $i_l \in [n_l],$ $l = 1, \cdots, k, k+p+1, \cdots, k+p+q$, with $k, p , q \in \mathbb{N}$.
Then,
$$\begin{aligned}
& \| \C \|_1 \\
& = \sum_{i_1=1}^{n_1} \cdots \sum_{i_k=1}^{n_k}\sum_{i_{k+p+1}=1}^{n_{k+p+1}}\cdots \sum_{i_{k+p+q}=1}^{n_{k+p+q}} \left|\sum_{i_{k+1}=1}^{n_{k+1}} \cdots \sum_{i_{k+p}=1}^{n_{k+p}}a_{i_1\cdots i_{k+p}}b_{i_{k+1}\cdots i_{k+p+q}}\right| \\
& \le \sum_{i_1=1}^{n_1} \cdots \sum_{i_k=1}^{n_k}\sum_{i_{k+p+1}=1}^{n_{k+p+1}}\cdots \sum_{i_{k+p+q}=1}^{n_{k+p+q}} \sum_{i_{k+1}=1}^{n_{k+1}} \cdots \sum_{i_{k+p}=1}^{n_{k+p}}\left|a_{i_1\cdots i_{k+p}}b_{i_{k+1}\cdots i_{k+p+q}}\right| \\
& \le \sum_{i_1=1}^{n_1} \cdots \sum_{i_k=1}^{n_k}\sum_{j_{k+p+1}=1}^{n_{k+p+1}}\cdots \sum_{j_{k+p+q}=1}^{n_{k+p+q}} \left(\sum_{i_{k+1}=1}^{n_{k+1}} \cdots \sum_{i_{k+p}=1}^{n_{k+p}}\left|a_{i_1\cdots i_{k+p}}\right|\right) \left(\sum_{j_{k+1}=1}^{n_{k+1}} \cdots \sum_{j_{k+p}=1}^{n_{k+p}}\left|b_{j_{k+1}\cdots  j_{k+p+q}}\right|\right)\\
& \le \left(\sum_{i_1=1}^{n_1} \cdots \sum_{i_{k+p}=1}^{n_{k+p}}\left|a_{i_1\cdots i_{k+p}}\right|\right) \left(\sum_{j_{k+1}=1}^{n_{k+1}} \cdots \sum_{j_{k+p+q}=1}^{n_{k+p+q}}\left|b_{j_{k+1}\cdots  j_{k+p+q}}\right|\right)\\
& = \|\A\|_1 \|\B\|_1.
\end{aligned}
$$
This shows that the 1-norm is a tensor norm.

Let $\A, \B$ and $\C$ be as defined above.  Denote the Frobenius norm as $\|\cdot\|_2$, as it is just the 2-norm.
Then
$$\|\A\|^2_2 = \sum_{i_1=1}^{n_1} \cdots \sum_{i_{k+p}=1}^{n_{k+p}}a_{i_1\cdots i_{k+p}}^2,$$
$$\|\B\|^2_2 = \sum_{i_{k+1}=1}^{n_{k+1}} \cdots \sum_{i_{k+p+q}=1}^{n_{k+p+q}}b_{i_{k+1}\cdots i_{k+p+q}}^2,$$
and
$$\| \C \|_2^2
= \sum_{i_1=1}^{n_1} \cdots \sum_{i_k=1}^{n_k}\sum_{i_{k+p+1}=1}^{n_{k+p+1}}\cdots \sum_{i_{k+p+q}=1}^{n_{k+p+q}} \left(\sum_{i_{k+1}=1}^{n_{k+1}} \cdots \sum_{i_{k+p}=1}^{n_{k+p}}a_{i_1\cdots i_{k+p}}b_{i_{k+1}\cdots i_{k+p+q}}\right)^2.$$
These show that
$$\|\C\|_2^2 \le \|\A\|_2^2 \|\B \|_2^2,$$
i.e.,
$$\|\C\|_2 \le \|\A\|_2 \|\B \|_2.$$
Hence, the Frobenius norm is also a tensor norm.
\qed

\section{The Cubic Power of a Third Order Tensor}

Consider the third order space $\Re^{d_1 \times d_2 \times d_3}$, where $d_1, d_2, d_3 \in \mathbb{N}$.    Such a third order tensor is the model of a general higher order tensor \cite{Hu15, YZ16}.

\begin{definition}  \label{d6.1}
Suppose that $\A = (a_{ijk}) \in \Re^{d_1 \times d_2 \times d_3}$.   Let $F :  \Re^{d_1 \times d_2 \times d_3} \to \Re^{d_1 \times d_2 \times d_3}$ be defined as follows.  We have $\A^3 \equiv F(\A) \equiv \T = (t_{ijk}) \in \Re^{d_1 \times d_2 \times d_3}$, where
\begin{equation} \label{eeee1}
t_{ijk} = \sum_{s=1}^{d_1} \sum_{p=1}^{d_2} \sum_{q=1}^{d_3} a_{ipq}a_{sjq}a_{spk},
\end{equation}
for $i \in [d_1], j \in [d_2]$ and $k \in [d_3]$.   In particular, we have
$$\A^{3^{m+1}} = F(\A^{3^m}),$$
for $m \in \mathbb{N}$.     If $\A^{3^m} = \O$ for some $m \in \mathbb{N}$, then $\A$ is said to be nilpotent.  If $\A^3 = \A$, then $\A$ is said to be idempotent.
\end{definition}

We see that $\A^3$ is uniquely defined for $\A \in \Re^{d_1 \times d_2 \times d_3}$.

\begin{example} \label{ex1}
Suppose that $d_1=d_2=d_3=d$.   We say that $\A = (a_{ijk}) \in \Re^{d \times d \times d}$ is diagonal if $a_{ijk} = 0$ as long as $i, j$ and $k$ are not all equal.  Assume that $a_{iii} = \alpha_i$ for $i \in [d]$.  Then $\A^3 \equiv \T = (t_{ijk}) \in \Re^{d \times d \times d}$ is also diagonal with $t_{iii} = \alpha^3_i$ for $i \in [d]$.
\end{example}

Thus, $\A^3$ preserves the diagonal property.   The matrix power preserves nonnegativity and symmetry.  We see $\A^3$ also preserves these.

\begin{proposition}
Suppose that $\A = (a_{ijk}) \in \Re^{d_1 \times d_2 \times d_3}$.    If $\A$ is nonnegative, i.e., $a_{ijk} \ge 0$ for $i \in [d_1], j \in [d_2], k \in [d_3]$, then $\A^3$ is also nonnegative.  If $d_1 = d_2 = d_3 =d$ and $\A$ is symmetric, i.e., $a_{ijk}$ is invariant under any index permutation, then $\A^3$ is also symmetric.
\end{proposition}
{\bf Proof} It is directly from the definition that $\A^3$ preserves nonnegativity.   Assume now $d_1 = d_2 = d_3 = d$ and $\A$ is symmetric. Denote $\T = (t_{ijk}) \equiv \A^3$.   For $i, j, k \in [d]$, we have
$$t_{jik} = \sum_{s, p, q =1}^d a_{jpq}a_{siq}a_{spk} = \sum_{s, p, q =1}^d a_{pjq}a_{isq}a_{psk} = t_{ijk}.$$
Then $\T = \A^3$ is also symmetric.
\qed

\begin{example} \label{ex2}
We consider the tensor $\A = (a_{ijk}) \in \Re^{2 \times 2 \times 2}$ with $a_{111} = -1$, $a_{122} = a_{212} = a_{221} = 1$, and $a_{ijk} = 0$ otherwise.
    Then $\A$ is symmetric.   We find that $\A^3 = \O$.   Thus, $\A$ is nilpotent.
  By computation, we find that the spectral norm of $\mathcal{A}$ is $1$.
  \end{example}

For a matrix $A \in \Re^{n \times n}$, if $A^k \not = 0$ for all positive integers $k \le n$, then $A$ is not nilpotent.   Is this also true for third order tensors?  That is, is there a number $L = L(d_1, d_2, d_3)$ such that if $\A^{3^m} \not = \O$ for all $m$ satisfying $3^m \le L$, then $\A \in \Re^{d_1 \times d_2 \times d_3}$ is not nilpotent?

Clearly, $\O$ is idempotent.   We now present some third order nontrivial idempotent tensors.

\begin{proposition} \label{p6.5}
Let $\x \in \Re^{d_1}$, $\y \in \Re^{d_2}$, $\z \in \Re^{d_3}$, such that $\| \x \|_2 = \|\y\|_2 = \|\z\|_2 = 1$.
Then $\x \circ \y \circ \z \in \Re^{d_1 \times d_2 \times d_3}$ is idempotent.
\end{proposition}
{\bf Proof} Checking by definition, we find that
$$\left(\x \circ \y \circ \z\right)^3 = \x \circ \y \circ \z.$$
\qed

Idempotent matrices have applications in regression analysis and econometrics.   Do third order idempotent tensors have some applications?

Coming back to tensor norms, we have the following proposition.

\begin{proposition} \label{p6.6}
Suppose $\A \in \Re^{d_1 \times d_2 \times d_3}$. If $|||\A||| < 1$ for a tensor norm $|||\cdot |||$, then
\begin{equation}  \label{ae1}
\lim_{m \to \infty} \A^{3^m} = \O.
\end{equation}
\end{proposition}
{\bf Proof} For a tensor norm $|||\cdot|||$, we have
\begin{equation}  \label{ae2}
|||\A^{3^m}||| \le |||\A|||^{3^m}.
\end{equation}
The conclusion follows. \qed

From Theorem \ref{t3}, $|||\cdot|||$ can be either the $1$-norm, or the Frobenius norm, or the nuclear norm.   By (\ref{ae2}), the convergence of (\ref{ae1}) should be very fast.

For a matrix $A \in \Re^{n \times n}$,
$$ \lim_{k \to \infty} A^k = 0$$
if and only if its spectral radius $\rho(A) < 1$.   See Theorem 5.6.12 of \cite{HJ12}.   In the next section, we will show that a similar result holds for third order tensors.

Another question is on the inverse operation of the cubic power.   For any $\A \in \Re^{d_1 \times d_2 \times d_3}$, does there exist $\B \in \Re^{d_1 \times d_2 \times d_3}$ such that $\A = \B^3$?   If such a $\B$ exists, is it unique?   If such a $\B$ exists and is unique, then we may call it the cubic root of $\A$ and denote
$$\B = \left(\A\right)^{1 \over 3}.$$

\subsection{The Power of an Odd Order Tensor}

We may extend Definition \ref{d6.1} to higher odd order tensors.   Let us consider the fifth order space $\Re^{d_1 \times d_2 \times d_3 \times d_4 \times d_5}$, where $d_1, d_2, d_3, d_4, d_5 \in \mathbb{N}$.

\begin{definition}  \label{d6.2}
Suppose that $\A = (a_{i_1i_2i_3i_4i_5}) \in \Re^{d_1 \times d_2 \times d_3 \times d_4 \times d_5}$.   Let $F :  \Re^{d_1 \times d_2 \times d_3 \times d_4 \times d_5} \to \Re^{d_1 \times d_2 \times d_3 \times d_4 \times d_5}$ be defined as follows.  We have $\A^5 \equiv F(\A) \equiv \T = (t_{i_1i_2i_3i_4i_5}) \in \Re^{d_1 \times d_2 \times d_3 \times d_4 \times d_5}$, where
\begin{equation} \label{eeee2}
t_{i_1i_2i_3i_4i_5} = \sum_{j_1, k_1=1}^{d_1} \sum_{j_2, k_2=1}^{d_2} \sum_{j_3, k_3=1}^{d_3} \sum_{j_4, k_4=1}^{d_4} \sum_{j_5, k_5=1}^{d_5} a_{i_1j_2j_3j_4j_5}a_{j_1i_2j_3j_4j_5}a_{j_1j_2i_3k_4k_5}a_{k_1k_2k_3i_4k_5}a_{k_1k_2k_3k_4i_5},
\end{equation}
for $i_l \in [d_l], l = 1, 2, 3, 4, 5$.   In particular, we have
$$\A^{5^{m+1}} = F(\A^{5^m}),$$
for $m \in \mathbb{N}$.     If $\A^{5^m} = \O$ for some $m \in \mathbb{N}$, then $\A$ is said to be nilpotent.   If $\A^5 = \A$, then $\A$ is said to be idempotent.
\end{definition}

We may still conduct our discussion.    However, such a definition is not unique now.   For example, we may replace (\ref{eeee2}) by
\begin{equation} \label{eeee3}
t_{i_1i_2i_3i_4i_5} = \sum_{j_1, k_1=1}^{d_1} \sum_{j_2, k_2=1}^{d_2} \sum_{j_3, k_3=1}^{d_3} \sum_{j_4, k_4=1}^{d_4} \sum_{j_5, k_5=1}^{d_5} a_{i_1j_2j_3j_4j_5}a_{k_1i_2j_3j_4j_5}a_{j_1j_2i_3k_4k_5}a_{j_1k_2k_3i_4k_5}a_{k_1k_2k_3k_4i_5},
\end{equation}

We may extend our discussion to all higher odd order tensors in this way.

\section{Gelfand Limit}

In matrix analysis, there is a well-known Gelfand formula (1941):
$$\rho(A)  = \lim_{k\to\infty} |||A^k|||^{1 \over k},$$
for any matrix norm $|||\cdot |||$.    See 5.6.14 of \cite{HJ12}.

Does the following limit
$$\lim_{m\to\infty} |||\A^{3^m}|||^{1 \over 3^m}$$
always exist for any $\A \in  \Re^{d_1 \times d_2 \times d_3}$ and any tensor norm  $|||\cdot |||$?   We now present the following theorem.

\begin{theorem} \label{t7.2}
Let $\A \in \Re^{d_1 \times d_2 \times d_3}$.   Then the limit
\begin{equation} \label{ee3}
\lim_{m \to \infty} \|\A^{3^m}\|^{1 \over 3^m}
\end{equation}
exists for any norm $\|\cdot \|$, and takes the same value for different norms.  If $\|\cdot \|$ is a tensor norm, then the sequence $\{ \|\A^{3^m}\|^{1 \over 3^m} \}$ is non-increasing.
\end{theorem}
{\bf Proof}  Let $\|\cdot\|$ be a tensor norm.    Then
$$\|\A^{3^{m+1}}\| \le \|\A^{3^m}\|^3,$$
i.e.,
$$\|\A^{3^{m+1}}\|^{1 \over 3^{m+1}} \le \|\A^{3^m}\|^{1 \over 3^m}.$$
This shows that the sequence $\{ \|\A^{3^m}\|^{1 \over 3^m} \}$ is non-increasing.   Since this sequence is nonnegative, it has a limit.

Now assume that $\|\cdot\|'$ is another norm of $\Re^{d_1 \times d_2 \times d_3}$, which is not necessary to be a tensor norm.  By the norm equivalence property in the finite dimensional space, there are positive constants $c_1$ and $c_2$ such that for any $\B \in \Re^{d_1 \times d_2 \times d_3}$,
$$c_1\|\B\| \le \|\B\|' \le c_2 \|\B\|.$$
Let $\B = \A^{3^m}$.    Then
$$c_1^{1\over 3^m} \|\A^{3^m}\|^{1\over 3^m} \le \left(\|\A^{3^m}\|'\right)^{1\over 3^m} \le c_2^{1 \over 3^m} \|\A^{3^m}\|^{1\over 3^m}.$$

Let $m \to \infty$.  Then we have
$$\lim_{m \to \infty} \|\A^{3^m}\|^{1\over 3^m} \le \liminf_{m \to \infty} \left(\|\A^{3^m}\|'\right)^{1\over 3^m} \le \limsup_{m \to \infty} \left(\|\A^{3^m}\|'\right)^{1\over 3^m} \le \lim_{m \to \infty} \|\A^{3^m}\|^{1\over 3^m}.$$
This shows that $\lim_{m \to \infty} \left(\|\A^{3^m}\|'\right)^{1\over 3^m}$ also exists and is equal to $\lim_{m \to \infty} {\|\A^{3^m}\|}^{1\over 3^m}$.

 Therefore, this limit exists for all the other norms with the same value.   This proves the theorem.
\qed

\begin{definition}
Let $\A \in \Re^{d_1 \times d_2 \times d_3}$.   We call $$\lim_{m \to \infty} \|\A^{3^m}\|^{1 \over 3^m}$$
the Gelfand limit of $\A$, and denote it as $\rho(\A)$.
\end{definition}

We check some properties of $\rho(\A)$.

\begin{theorem}
Let $\A \in \Re^{d_1 \times d_2 \times d_3}$ and $\alpha \in \Re$.    Then we have
$$\rho(\alpha\A) = |\alpha|\rho(\A)$$
and
$$\rho(\A^3) = \rho(\A)^3.$$

If $\|\cdot\|$ is a tensor norm, then
$$\rho(\A) \le \|\A\|.$$

The Gelfand limit is zero if $\A$ is nilpotent.   The Gelfand limit is one or zero if $\A$ is idempotent.
\end{theorem}
{\bf Proof}     By definition, we have
$$\rho(\alpha\A) = \lim_{m\to\infty} \left\|(\alpha\A)^{3^m}\right\|^{1 \over 3^m} = |\alpha|\lim_{m\to\infty} ||\A^{3^m}||^{1 \over 3^m} = |\alpha|\rho(\A)$$
and
$$\rho(\A^3) = \lim_{m\to\infty} \left\|\A^{3^{m+1}}\right\|^{1 \over 3^m} = \lim_{m\to\infty} \left(||\A^{3^{m+1}}||^{1 \over 3^{m+1}}\right)^3 = \rho(\A)^3.$$

The second property follows from Theorem \ref{t7.2}.

If $\A$ is nilpotent, then there is $M \in \mathbb{N}$ such that
$$\A^{3^m} \equiv \O$$
for all $m \ge M$.   Then $\rho(\A) = 0$.

If $\A$ is idempotent, then
$$\rho(\A) = \rho(\A^3) = \rho(\A)^3.$$
Then $\rho(\A) = 1$ or $0$.
 \qed

If $\rho(\A) =0$, is $\A$ nilpotent?  If $\A$ is idempotent and $\rho(\A) = 0$, do we have $\A = \O$?

\begin{theorem}
Suppose $\A \in \Re^{d_1 \times d_2 \times d_3}$. Then
\begin{equation}  \label{aae1}
\lim_{m \to \infty} \A^{3^m} = \O,
\end{equation}
if and only if  $\rho(\A) < 1$.
\end{theorem}
{\bf Proof} Suppose that $\rho(\A) < 1$.  Let $\|\cdot\|$ be a tensor norm.   Then there are $M \in \mathbb{N}$ such that
$$\|\A^{3^M}\|^{1 \over 3^M} < 1.$$
This implies that
$$\|\A^{3^M}\| < 1.$$
By Proposition \ref{p6.6}, we have
$$\lim_{m \to \infty} \left(\A^{3^M}\right)^{3m} = \O.$$
Then (\ref{aae1}) follows.

On the other hand, suppose that $\rho(\A) \ge 1$.   Let $\|\cdot\|$ be a tensor norm.  By Theorem \ref{t7.2}, $\left\{ \|\A^{3^m} \|^{1 \over 3^m} \right\}$ is non-increasing and tends to $\rho(\A)$.  Hence, for any $m \in \mathbb{N}$,
$$\|\A^{3^m} \|^{1 \over 3^m} \ge \rho(\A) \ge 1,$$
which implies that
$$\|\A^{3^m} \| \ge 1.$$
Then (\ref{aae1}) cannot hold.
\qed

  \begin{example}  We generate a general tensor
  $\mathcal{A} = (a_{ijk}) \in\Re^{4\times 2\times2}$ randomly:
  \begin{equation*}
    \left(
       \begin{array}{cc|cc}
  a_{111} &   a_{121} &  a_{112} &  a_{122} \\
   a_{211} &   a_{221} &  a_{212} & a_{222} \\
   s_{311} &  a_{321} &   a_{312} &  a_{322} \\
   s_{411} &  a_{421} &   a_{412} &  a_{422} \\
       \end{array}
     \right).
  \end{equation*}
    \begin{equation*}
    =\left(
       \begin{array}{ccc|ccc}
  -0.512159 &  -0.507535 &  -1.383216 &   0.203856 &  -0.578312 &   1.921669 \\
   0.906334 &  -0.258462 &  -0.982083 &   0.736707 &   0.608575 &  -1.063641 \\
  -0.731184 &   0.525138 &  -1.347676 &  -0.782006 &   0.568222 &  -0.214013 \\
  -0.086664 &  -0.736508 &   0.474856 &   0.345770 &   0.194509 &   0.006420 \\
       \end{array}
     \right).
  \end{equation*}
 The sequences of $\|\mathcal{A}^{3^m}\|_1$, $\|\mathcal{A}^{3^m}\|^{1/3^m}_2$ and $\|\mathcal{A}^{3^m}\|_\infty^{1/3^m}$  are: \\
  \begin{center}
  \begin{tabular}{c|lll}
    $m$ & $\|\mathcal{A}^{3^m}\|_1^{1/3^m}$  & $\|\mathcal{A}^{3^m}\|_2^{1/3^m}$ & $\|\mathcal{A}^{3^m}\|_{\infty}^{1/3^m}$ \\
    \hline
0 & 15.6755 & 3.86508 & 1.92167 \\

1 & 4.0199 & 2.70142 & 2.31202 \\

2 & 2.82591 & 2.54596 & 2.45592 \\

3 & 2.61624 & 2.53718 & 2.50769 \\

4 & 2.56299 & 2.53712 & 2.52722 \\

5 & 2.54571 & 2.53712 & 2.53382 \\

6 & 2.53998 & 2.537118666456933 & 2.53602 \\

7 & 2.538072064165983 & 2.537118666456933 & 2.536751440470295 \\

8 & 2.537436425894090 & 2.537118666456933 & 2.536996251888380 \\

9 & 2.537224581847678 & 2.537118666456933 & 2.537077860944460 \\

10 & 2.537153971095906 & 2.537118666456933 & 2.537105064546520 \\

11 & 2.537130434615338 & 2.537118666456933 & 2.537114132478693 \\

12 & 2.537122589170336 & 2.537118666456933 & 2.537117155129952 \\

13 & 2.537119974027393 & 2.537118666456933 & 2.537118162681173 \\

14 & 2.537119102313678 & 2.537118666456933 & 2.537118498531668 \\

15 & 2.537118811742506 & 2.537118666456933 & 2.537118610481843 \\

16 & 2.537118714885456 & 2.537118666456933 & 2.537118647798569 \\

17 & 2.537118682599774 & 2.537118666456933 & 2.537118660237478 \\

18 & 2.537118671837880 & 2.537118666456933 & 2.537118664383781 \\

19 & 2.537118668250582 & 2.537118666456933 & 2.537118665765882 \\

20 & 2.537118667054816 & 2.537118666456933 & 2.537118666226583 \\

21 & 2.537118666656227 & 2.537118666456933 & 2.537118666380149 \\

22 & 2.537118666523364 & 2.537118666456933 & 2.537118666431338 \\

23 & 2.537118666479076 & 2.537118666456933 & 2.537118666448401 \\

24 & 2.537118666464314 & 2.537118666456933 & 2.537118666454089 \\

25 & 2.537118666459393 & 2.537118666456933 & 2.537118666455985 \\

26 & 2.537118666457753 & 2.537118666456933 & 2.537118666456617 \\

27 & 2.537118666457206 & 2.537118666456933 & 2.537118666456827 \\

28 & 2.537118666457024 & 2.537118666456933 & 2.537118666456898 \\

29 & 2.537118666456963 & 2.537118666456933 & 2.537118666456921 \\

30 & 2.537118666456943 & 2.537118666456933 & 2.537118666456929 \\

31 & 2.537118666456936 & 2.537118666456933 & 2.537118666456931 \\
  \end{tabular}
  \end{center}
  The three sequences converge to the same limit.   The first two sequences are non-increasing, as $\|\cdot\|_1$ and $\|\cdot\|_2$ norms are tensor norms.
  \end{example}

We see that $\rho(\A)$  clearly distinguishes tensor norms from the other norms.

In matrix analysis \cite{HJ12}, for any given matrix $A \in \Re^{n \times n}$ and $\epsilon > 0$, there is a matrix norm $|||\cdot |||$ such that
$$\rho(A) \le |||A||| \le \rho(A) + \epsilon.$$
Is this true for tensors?   That is, for any given $\A \in \Re^{d_1 \times d_2 \times d_3}$ and $\epsilon > 0$, is there a tensor norm $|||\cdot |||$ such that
$$\rho(\A) \le |||\A||| \le \rho(\A) + \epsilon?$$

There are various mysterious things surrounding the Gelfand limit of a third order tensor.   After unveiling such mysterious things, we should understand the spectral theory of tensors more.

\subsection{Higher Odd Order Tensors}

By applying the discussion in Subsection 6.1, we may extend the discussion of the Gelfand limit to higher order tensors.    The question is: if we use (\ref{eeee3}) instead of (\ref{eeee2}), does the Gelfand limit have the same value?  Preliminary computation shows that the Gelfand limit may have different values in this case.  Then, the definitions of higher odd order nilpotent and idempotent tensors may also vary.  However, in this case, the higher order extension of Proposition \ref{p6.5} is still valid. This means that some essential things may be invariant under (\ref{eeee2})
and (\ref{eeee3}).

%\section{Universal Vector Norm and Induced Tensor Norm}

%We now introduce the concept of universal vector norm.

%\begin{definition}
%Suppose that $\| \cdot \| : \R^n \to \Re_+$ for any $n \in \mathbb{N}$ such that for any fixed $n \in \mathbb{N}$, it is a norm on $\Re^n$; and if $m > n$, $m, n \in \mathbb{N}$, then for any $\x = (x_1, \cdots, x_n)^\top \in \Re^n$, we have
%$$\|\x \| = \| \y \|,$$
%where $\y = (x_1, \cdots, x_n, 0, \cdots, 0)^\top \in \Re^m$, i.e., $\x$ is extended to $\y$ by adding $n-m$ $0$'s.
%Then $\| \cdots \|$ is a universal vector norm.
%\end{definition}

\section{Final Remarks}

In this paper, we we extended the concept of matrix norm to tensor norm.   We showed that the $1$-norm, the Frobenius norm and the nuclear norm of tensors are tensor norms but the infinity norm and the spectral norm of tensors are not tensor norms.

We introduced the cubic power for a general third order tensor, and showed that the cubic power of a general third order tensor tends to zero as the power increases to infinity, if there is a tensor norm such that the tensor norm of that third order tensor is less than one.   Then we showed that a Gelfand formula holds for a general third order tensor.
We call such a limit the Gelfand limit of that third order tensor. The Gelfand limit is zero for all third order nilpotent tensors, and one or zero for all third order idempotent tensors. We showed that the Gelfand limit is less than or equal to any tensor norm value of that third order tensor, and the cubic power of that third order tensor tends to zero as the power increases to infinity if and only if the Gelfand limit is less than one.   Our results on the cubic power and the Gelfand limit can be extended to any higher odd order tensors.   As the limit in the matrix Gelfand formula is the spectral radius of the matrix, this may open the way for a new spectral theory for odd order tensors.

Our derivations are on the field of real numbers.  They can be extended to the field of complex numbers without difficulty.

\end{document}